\documentclass[11pt]{article}
\usepackage{amsmath,amssymb,amsthm}
\usepackage[margin=1in]{geometry}
\usepackage{hyperref}
\usepackage{tikz-cd}

\newtheorem{theorem}{Theorem}
\newtheorem{lemma}{Lemma}
\newtheorem{corollary}{Corollary}
\theoremstyle{definition}
\newtheorem{definition}{Definition}

\newcommand{\Z}{\mathbb{Z}}
\newcommand{\F}{\mathbb{F}}

\title{Unary maps are restrictions of linear maps and $k$-ary maps are restrictions of
homogeneous maps}
\author{Roman Bacik}
\date{November 29, 2026}

\begin{document}
\maketitle

\begin{abstract}
We study maps on a finite set $S$ under renaming the elements of $S$ and extending $S$ to a ring.
We show that every map $S\to S$ is the restriction of a linear map on some $\mathbb{Z}/m$, and
every map $S^k\to S$ is the restriction of a homogeneous polynomial map over some prime field.
Nothing simpler works in general: the number of separable terms and the degree cannot be
bounded, and the degree must grow with $k$.
\end{abstract}

\section{Introduction}

A cellular automaton updates all cells at once by a local rule $f\colon S^k\to S$ that reads the
states of $k$ neighbouring cells, and already two-neighbour automata can be
computation-universal~\cite{smith}. We look for the simplest functional description of an arbitrary rule $f$, while
permitting renaming and enlarging of the set of states: the states are embedded into a ring $R$
and $f$ becomes the restriction of a simple map on $R$. This changes neither the transitions
between the original states nor what the automaton computes.

Embedding into an arbitrary polynomial is trivial: for a prime $p\ge|S|$ every injective map
$j\colon S\to\F_p$ works, since every map $\F_p^k\to\F_p$ is a polynomial~\cite{lidl}. We therefore ask for embeddings into polynomials of a restricted form. Every unary map
embeds into a linear map over some $\Z/m$ (Theorem~\ref{thm:unary}) and into a monomial over some
prime field (Lemma~\ref{lem:power}). Theorem~\ref{thm:binary} extends the second statement to all
arities: every map $S^k\to S$ embeds into a homogeneous polynomial over some prime field.

This note continues~\cite{finlin} and~\cite{finbin}. Theorem~\ref{thm:unary} is the main result
of~\cite{finlin}. In~\cite{finbin} it is shown that for every $d$ some binary map is not represented
by a polynomial of degree $d$ over any commutative ring; Theorem~\ref{thm:optimal}(b) is from
there. New here are Theorem~\ref{thm:binary}, the bounds of
Theorem~\ref{thm:optimal}(a) and Corollary~\ref{cor:degree}, and Theorem~\ref{thm:arity}, which
shows that the degree must grow with $k$.

Polynomial functions on $\Z/m$ go back to Kempner~\cite{kempner}; Theorem~\ref{thm:unary} shows
that every map on a finite set is the restriction of a linear function on some $\Z/m$. In
complexity theory a Boolean function is represented by a polynomial with a fixed encoding of the
inputs and outputs. Over the reals the degree does not depend on which two numbers encode
\emph{true} and \emph{false}~\cite{nisanszegedy}. Modulo $m$, with the weaker requirement that
only zero and nonzero outputs are distinguished, the OR of $N$ variables has degree
$O(N^{1/r})$, where $r$ is the number of distinct prime factors of $m$, while modulo a prime $p$
it has degree $\lceil N/(p-1)\rceil$~\cite{bbr}. Here both the encoding $j$ and the ring are
free, but outputs are encoded exactly, and Theorem~\ref{thm:arity} shows that the conjunction of
$k$ variables still needs degree $k$ over every commutative ring. The number of separable terms
in Theorem~\ref{thm:optimal}(a) is a rank of the matrix $(j(f(x,y)))_{x,y}$; the rank of such
matrices bounds deterministic communication complexity~\cite{mehlhornschmidt}.

Let $S$ be a finite set and $R$ a commutative ring.

\begin{definition}
Let $k\ge 1$. An \emph{embedding} of a map $f\colon S^k\to S$ into a map $q\colon R^k\to R$ is an
injective map $j\colon S\to R$ with $j\circ f=q\circ j^k$, where
$j^k(x_1,\dots,x_k)=(j(x_1),\dots,j(x_k))$, that is, such that the square
\[
\begin{tikzcd}
S^k \arrow[r,"f"] \arrow[d,hook,"j^k"'] & S \arrow[d,hook,"j"] \\
R^k \arrow[r,"q"'] & R
\end{tikzcd}
\]
commutes. For $k=1$ and $k=2$ the map $f$ is \emph{unary} and \emph{binary}, and for binary $f$
we write $f(x,y)$ and $q(u,v)$.
\end{definition}

If $f$ has an embedding into $q$, we say that $f$ is \emph{represented} by $q$ over $R$.

\begin{definition}
A map $q\colon R\to R$ is \emph{linear} if $q(u)=c\,u$ for some $c\in R$. A map
$q\colon R^2\to R$ is \emph{$r$-separable}, $r\ge 0$, if $q(u,v)=\sum_{i=1}^{r}g_i(u)\,h_i(v)$
for some maps $g_i,h_i\colon R\to R$.
\end{definition}

Identifying $S$ with $j(S)\subseteq R$, an embedding of $f$ into $q$ makes $f$ the
\emph{restriction} of $q$ to $j(S)^k$.

\section{Results}

Throughout, $s=|S|$, $\Z/0=\Z$, and $[P]\in\{0,1\}$ is the truth value of a statement $P$; we write
$[P]_R$ when the arithmetic in $P$ is evaluated in a ring $R$.

\begin{theorem}\label{thm:unary}
Every map $f\colon S\to S$ is represented by a linear map over some $\Z/m$, $m\ge 1$.
\end{theorem}

\begin{theorem}\label{thm:binary}
Let $k\ge 1$. Every map $f\colon S^k\to S$ is represented by a homogeneous polynomial over some
$\F_p$, $p$ prime.
\end{theorem}

To show that Theorem~\ref{thm:binary} cannot be simplified we use the shift, equality and
inversion indicators. For $n\ge 2$ we identify $\Z/n$ with $\{0,\dots,n-1\}$ and let
\[
  S_n(x,y)=[x=y+1]_\Z,\qquad E_n(x,y)=[x=y],\qquad I_n(x,y)=[xy=1]_{\Z/n}\qquad(x,y\in\Z/n).
\]

\begin{theorem}\label{thm:optimal}
\begin{enumerate}
\item[(a)] Let $n\ge 2$ and $m\ge 0$. If $m=0$ or $m\ge n$, every binary map on an $n$-element
set is represented by an $n$-separable map over $\Z/m$. The maps $S_n$ and $E_n$ are not
represented by an $(n-2)$-separable map over $\Z/m$.
\item[(b)] For every prime $p$ and every commutative ring $R$, the map $I_p$ is not represented
over $R$ by a polynomial of total degree at most $p-1$.
\end{enumerate}
\end{theorem}

Part (b) is~\cite[Theorem~5]{finbin}; we include its proof for completeness.

\begin{corollary}\label{cor:degree}
Let $n\ge 2$ and $m\ge 0$. The maps $S_n$ and $E_n$ are not represented over $\Z/m$ by a
polynomial of total degree at most $n-2$.
\end{corollary}

For $k\ge 1$ let $A_k\colon(\Z/2)^k\to\Z/2$ be the conjunction $A_k(x)=x_1\cdots x_k$.

\begin{theorem}\label{thm:arity}
Let $k\ge 1$. For every commutative ring $R$, the map $A_k$ is not represented over $R$ by a
polynomial of total degree at most $k-1$. Over $\Z$, the map $A_k$ is represented by the
homogeneous polynomial $u_1\cdots u_k$.
\end{theorem}

\section{Proof of Theorem~\ref{thm:unary}}

A map $\varphi\colon S\to G$ into an abelian group is \emph{Sidon}~\cite{obryant} if
$\varphi(x)-\varphi(y)=\varphi(x')-\varphi(y')$ with $x\ne y$ and $x'\ne y'$ implies
$(x,y)=(x',y')$. A Sidon map is injective.

\begin{lemma}\label{lem:pow2}
Let $i,i',l,l'\ge 0$ with $i\ne i'$ and $l\ne l'$. If $2^i-2^{i'}=2^l-2^{l'}$, then $i=l$ and
$i'=l'$.
\end{lemma}

\begin{proof}
Both sides have the same sign, so we may assume $i>i'$ and $l>l'$. Then
$2^{i'}(2^{i-i'}-1)=2^{l'}(2^{l-l'}-1)$ with odd second factors, so $i'=l'$, and then $i=l$.
\end{proof}

\begin{lemma}\label{lem:sidon}
For every $f\colon S\to S$ and $d_0\ge 0$ there are $n\ge 1$, $d\ge d_0$ and a Sidon map
$\varphi\colon S\to\Z/n$ with $\varphi(f(x))=d\,\varphi(x)$ for all $x\in S$.
\end{lemma}

This refines~\cite[Lemma~2.12]{finlin}, which uses the weights $i+1$ in place of $2^i$ and
obtains an injective map; the weights $2^i$ make it Sidon.

\begin{proof}
If $s\le 1$, take $n=1$, any $d\ge d_0$ and $\varphi=0$. Let $s\ge 2$, identify $S=\{0,\dots,s-1\}$,
and let $A$ be the $s\times s$ matrix $A_{ik}=[f(i)=k]$, so that $(Az)_i=z_{f(i)}$. Put
$\beta_i=2^i$, $\chi(t)=\det(tI-A)\in\Z[t]$ and $w(t)=\operatorname{adj}(tI-A)\,\beta\in\Z[t]^s$.
From $(tI-A)\operatorname{adj}(tI-A)=\chi(t)I$ we get $Aw=tw-\chi(t)\beta$, that is,
\[
  w_{f(i)}(t)=t\,w_i(t)-\chi(t)\,\beta_i .
\]
The polynomial $\chi$ is monic of degree $s$. The diagonal entries of
$\operatorname{adj}(tI-A)$ are monic of degree $s-1$ and the other entries have degree at most
$s-2$~\cite[Lemma~2.7]{finlin},~\cite{horn}, so $w_i$ has degree at most $s-1$ and coefficient $\beta_i$ at $t^{s-1}$. For $i\ne i'$,
$l\ne l'$ and $(i,i')\ne(l,l')$, the polynomial $w_i-w_{i'}-w_l+w_{l'}$ has coefficient
$\beta_i-\beta_{i'}-\beta_l+\beta_{l'}\ne 0$ at $t^{s-1}$ by Lemma~\ref{lem:pow2}. Hence for
every large integer $t_0$ we have $0<|w_i(t_0)-w_{i'}(t_0)-w_l(t_0)+w_{l'}(t_0)|<\chi(t_0)$ for
all such quadruples. Fix such a $t_0\ge\max(d_0,2)$ and put $n=\chi(t_0)$, $d=t_0$ and
$\varphi(i)=w_i(t_0)\bmod n$. The displayed identity gives $\varphi(f(i))=d\,\varphi(i)$, and
$\varphi$ is Sidon because $n$ divides no nonzero integer of absolute value less than $n$.
\end{proof}

\begin{proof}[Proof of Theorem~\ref{thm:unary}]
Take $n$, $d$ and $\varphi$ from Lemma~\ref{lem:sidon} and put $j=\varphi$. Then $j$ is
injective and $j\circ f=q\circ j$ for $q(u)=d\,u$, so $j$ is an embedding of $f$ into the linear
map $q$ over $\Z/n$.
\end{proof}

\section{Proof of Theorem~\ref{thm:binary}}

\begin{lemma}\label{lem:power}
For every $f\colon S\to S$ and $d_0\ge 0$ there are a prime $p$, an integer $d\ge d_0$ and an
embedding $j\colon S\to\F_p^\times$ of $f$ into $q(u)=u^d$ over $\F_p$ such that the ratios
$j(x)/j(y)$, $x\ne y$, are pairwise distinct.
\end{lemma}

\begin{proof}
Take $n$, $d$ and $\varphi$ from Lemma~\ref{lem:sidon}, a prime $p\equiv 1\pmod n$~\cite{borevich},
and $h\in\F_p^\times$ of order $n$~\cite{lang}. Put $j(x)=h^{\varphi(x)}$. Then
$j(f(x))=h^{d\,\varphi(x)}=j(x)^d$. Since $k\mapsto h^k$ is injective on $\Z/n$, so is
$j$, and the ratios $j(x)/j(y)=h^{\varphi(x)-\varphi(y)}$, $x\ne y$, are pairwise
distinct because $\varphi$ is Sidon.
\end{proof}

\begin{proof}[Proof of Theorem~\ref{thm:binary}]
Apply Lemma~\ref{lem:power} to the diagonal $a\mapsto f(a,\dots,a)$ with $d_0=s^k$. This gives a
prime $p$, an integer $d\ge s^k$ and an embedding $j\colon S\to\F_p^\times$ with
$j(f(a,\dots,a))=j(a)^d$ and pairwise distinct ratios $j(a)/j(b)$, $a\ne b$. For
$x=(x_1,\dots,x_k)\in S^k$ put
\[
  P(x)=\Bigl(\frac{j(x_1)}{j(x_1)},\dots,\frac{j(x_k)}{j(x_1)}\Bigr)\in\F_p^k,\qquad
  V(x)=j(f(x))\,j(x_1)^{-d}.
\]
If $P(x)=P(x')$, then $V(x)=V(x')$. Indeed, $j(x_i)j(x'_1)=j(x'_i)j(x_1)$ for all $i$. If
$x_1=x'_1$, then $x=x'$. Otherwise $x_i\ne x'_i$ for every $i$, since $x_i=x'_i$ would give $j(x'_1)=j(x_1)$, and the distinct ratios give
$x_i=x_1$ and $x'_i=x'_1$, so $x$ and $x'$ are constant; then $V(x)=V(x')=1$ by the choice of $j$.

Let $T=P(S^k)$, so $|T|\le s^k\le d$, and let $V_a=V(x)$ for any $x$ with $P(x)=a\in T$. For
$a\ne b$ in $T$ choose an index $c=c(a,b)$ with $a_c\ne b_c$, and put
\[
  L_{a,b}(u)=\frac{u_c-b_c\,u_1}{a_c-b_c},\qquad B_a=\prod_{b\in T\setminus\{a\}}L_{a,b},\qquad
  q=\sum_{a\in T}V_a\,u_1^{\,d-|T|+1}B_a .
\]
Each $L_{a,b}$ is a linear form, so $q$ is homogeneous of degree $d$. At $u=j^k(x)$ we have
$L_{a,b}(u)=j(x_1)\,(P(x)_c-b_c)/(a_c-b_c)$, which is $j(x_1)$ if $P(x)=a$ and $0$ if $b=P(x)$.
Hence $B_a(j^k(x))=j(x_1)^{|T|-1}\,[P(x)=a]$ and
$q(j^k(x))=V(x)\,j(x_1)^d=j(f(x))$. So $j$ is an embedding of $f$ into $q$.
\end{proof}

\section{Proof of Theorem~\ref{thm:optimal}}

We use a rank criterion.

\begin{lemma}\label{lem:lift}
Let $m,r\ge 0$, let $P$ be an integer matrix, let $c\in\Z/m$ be nonzero, and suppose $cP=EH$
over $\Z/m$ for matrices $E$ and $H$ over $\Z/m$ with $r$ columns and $r$ rows respectively. Then
$P$ has rank at most $r$ over $\F_p$ for some prime $p$.
\end{lemma}

\begin{proof}
Lift $E$, $H$ and $c$ to integer matrices $E'$, $H'$ and an integer $c'$. Every entry of
$E'H'-c'P$ is divisible by $m$, and $m\nmid c'$, so there are a prime $p$ and $k\ge 1$ with
$p^k\mid m$ and $p^k\nmid c'$. Write $c'=p^ec''$ with $p\nmid c''$; then $e<k$, so
$p^{e+1}\mid m$ and $E'H'=p^eB$ with $B=c''P+pD$ for an integer matrix $D$. The
$(r+1)\times(r+1)$ minors of $E'H'$ vanish, since $E'H'$ factors through $\mathbb{Q}^r$; they are
$p^{e(r+1)}$ times the minors of $B$. So $B\equiv c''P$ has rank at most $r$ over $\F_p$, and
$c''$ is a unit modulo $p$.
\end{proof}

\begin{lemma}\label{lem:criterion}
Let $f\colon X\times Y\to Z$, $r\ge 0$, $x_0\in X$, finitely many $x_i\in X$ and $y_l\in Y$,
and $u\ne w$ in $Z$ with $f(x_0,y_l)=u$ and $f(x_i,y_l)\in\{u,w\}$ for all $i,l$, and suppose that
the matrix $P_{il}=[f(x_i,y_l)=w]$ has rank greater than $r$ over $\F_p$ for every
prime $p$. Then for no $m\ge 0$ are there an injective map $j\colon Z\to\Z/m$ and maps
$G_k\colon X\to\Z/m$ and $H_0,H_k\colon Y\to\Z/m$ with
$j(f(x,y))=H_0(y)+\sum_{k=1}^{r}G_k(x)H_k(y)$ for all $x,y$.
\end{lemma}

\begin{proof}
Put $c=j(u)-j(w)\ne 0$, $E_{ik}=G_k(x_0)-G_k(x_i)$ and $H'_{kl}=H_k(y_l)$. The terms
$H_0(y_l)$ cancel in $j(f(x_0,y_l))-j(f(x_i,y_l))$, so $(EH')_{il}=c\,P_{il}$. By
Lemma~\ref{lem:lift} the matrix $P$ has rank at most $r$ over some $\F_p$, a contradiction.
\end{proof}

\begin{proof}[Proof of (a)]
Let $f$ be a binary map on an $n$-element set $S$. Since $m=0$ or $m\ge n$, there is an
injective map $j\colon S\to\Z/m$. For $z\in S$ put
$g_z(u)=[u=j(z)]$, $h_z(j(y))=j(f(z,y))$ for $y\in S$, and $h_z(v)$ arbitrary for
$v\notin j(S)$. Then $j$ is an embedding of $f$ into
$q(u,v)=\sum_{z\in S}g_z(u)\,h_z(v)$.

If $j$ is an embedding of $S_n$ into $\sum_{k=1}^{n-2}g_k(u)h_k(v)$, then
Lemma~\ref{lem:criterion} applies with $G_k(x)=g_k(j(x))$, $H_k(y)=h_k(j(y))$, $H_0=0$,
$r=n-2$, $x_0=0$, $x_i=i+1$ and $y_l=l$ for $0\le i,l\le n-2$, $u=0$ and $w=1$: then $P$ is the
$(n-1)\times(n-1)$ identity matrix, of rank $n-1>r$ over every $\F_p$, a contradiction. For $E_n$
use $y_l=l+1$ instead.
\end{proof}

\begin{proof}[Proof of Corollary~\ref{cor:degree}]
A polynomial of total degree at most $n-2$ has the form $\sum_{k=0}^{n-2}u^kq_k(v)$. If $j$ is an
embedding of $f=S_n$ or $f=E_n$ into it, then $j(f(x,y))=H_0(y)+\sum_{k=1}^{n-2}G_k(x)H_k(y)$ with
$G_k(x)=j(x)^k$, $H_k(y)=q_k(j(y))$ and $H_0(y)=q_0(j(y))$. This is excluded by
Lemma~\ref{lem:criterion} with the data of the proof of (a).
\end{proof}

For (b) let $\sigma_k(u_0,u_1)=\sum_{i=0}^{k-1}u_1^iu_0^{k-1-i}$, so that
$u_1^k-u_0^k=(u_1-u_0)\sigma_k(u_0,u_1)$, and let $q_{k_1k_2}$ be the coefficient of
$u^{k_1}v^{k_2}$ in $q$.

\begin{lemma}\label{lem:mixed}
If $q\in R[u,v]$ has total degree at most $N+1$, then for all $u_0,u_1,v_0,v_1\in R$
\[
  q(u_1,v_1)-q(u_0,v_1)-q(u_1,v_0)+q(u_0,v_0)=(u_1-u_0)(v_1-v_0)\sum_{k_1,k_2=1}^{N}
  q_{k_1k_2}\,\sigma_{k_1}(u_0,u_1)\,\sigma_{k_2}(v_0,v_1).
\]
\end{lemma}

\begin{proof}
The monomial $u^{k_1}v^{k_2}$ contributes $(u_1^{k_1}-u_0^{k_1})(v_1^{k_2}-v_0^{k_2})$, which
vanishes unless $k_1,k_2\ge 1$, and then $k_1,k_2\le N$.
\end{proof}

\begin{proof}[Proof of (b)]
Suppose $j$ is an embedding of $I_p$ into $q\in R[u,v]$ of total degree at most $p-1$. Put
$N=p-2$, $\delta=j(1)-j(0)$ and $U=(\Z/p)^\times$. For
$x,y\in U$, Lemma~\ref{lem:mixed} with $u_0=v_0=j(0)$, $u_1=j(x)$ and $v_1=j(y)$ gives, since
$I_p(0,y)=I_p(x,0)=I_p(0,0)=0$,
\begin{gather*}
  (j(x)-j(0))(j(y)-j(0))\,e(x,y)=j(I_p(x,y))-j(0)=\delta\,[xy=1]_{\Z/p},\\
  e(x,y)=\sum_{k_1,k_2=1}^{N}q_{k_1k_2}\,\sigma_{k_1}(j(0),j(x))\,\sigma_{k_2}(j(0),j(y)).
\end{gather*}
In matrix form, indexed by $U$, this reads $D\,\Phi Q\Phi^{\mathsf T}D=\delta F$, where
$D=\operatorname{diag}(j(x)-j(0))$, $\Phi_{x,k}=\sigma_k(j(0),j(x))$ for $1\le k\le N$,
$Q=(q_{k_1k_2})$, and $F$ is the permutation matrix of $x\mapsto x^{-1}$. Since
$\Phi Q\Phi^{\mathsf T}$ factors through $R^N$ and $N<p-1=|U|$, its determinant vanishes, while
$\det F=\pm1$. Taking determinants gives $\delta^{p-1}=0$. The entry $x=y=1$ gives
$\delta(1-\delta\,e(1,1))=0$; as $\delta\,e(1,1)$ is nilpotent, $1-\delta\,e(1,1)$ is a unit, so
$\delta=0$, contradicting injectivity.
\end{proof}

\section{Proof of Theorem~\ref{thm:arity}}

For $a,b\in R$ and $t\subseteq\{1,\dots,k\}$ let $v^t\in R^k$ have $v^t_l=b$ for $l\in t$ and
$v^t_l=a$ otherwise, and for $F\colon R^k\to R$ let
\[
  \Delta(F)=\sum_{t\subseteq\{1,\dots,k\}}(-1)^{k-|t|}\,F(v^t).
\]

\begin{lemma}\label{lem:kdiff}
If $q\in R[u_1,\dots,u_k]$ has coefficients $q_m$, then
$\Delta(q)=\sum_m q_m\prod_{l=1}^{k}\bigl(b^{m_l}-a^{m_l}\bigr)$. In particular $\Delta(q)=0$ if $q$
has total degree at most $k-1$.
\end{lemma}

\begin{proof}
For a monomial, $\prod_l(b^{m_l}-a^{m_l})$ expands to $\sum_t(-1)^{k-|t|}\prod_l(v^t_l)^{m_l}$. If
the total degree is at most $k-1$, every monomial misses some variable $u_l$, and its factor
$b^0-a^0$ is $0$.
\end{proof}

\begin{proof}[Proof of Theorem~\ref{thm:arity}]
Let $j$ be an embedding of $A_k$ into $q$ of total degree at most $k-1$, and put $a=j(0)$,
$b=j(1)$. Then $q(v^t)=j(A_k(x))$ with $x_l=[l\in t]$, so $q(v^t)=a+(b-a)\,[t=\{1,\dots,k\}]$ and
\[
  \Delta(q)=a\sum_t(-1)^{k-|t|}+(b-a)=a\,(1-1)^k+(b-a)=b-a .
\]
By Lemma~\ref{lem:kdiff}, $\Delta(q)=0$, so $j(1)=j(0)$, contradicting injectivity. Over $\Z$ take
$j(0)=0$, $j(1)=1$ and $q=u_1\cdots u_k$.
\end{proof}

\section*{Acknowledgements}
All results are formalized in Lean~4~\cite{lean4} with Mathlib~\cite{mathlib}; the formalization
is available at~\cite{hombin}.
Claude Opus~\cite{claude} and Aristotle~\cite{aristotle} were used for the Lean formalization. 
The mathematics is the author's,
and all proofs were checked by the author and formally verified.

\end{document}